# The strict complementary slackness condition in linear fractional programming

Mahmood Mehdiloozad

Kaoru Tone

Mohammad Bagher Ahmadi

March 2016

(✉) **Mahmood Mehdiloozad**
Department of Mathematics, College of Sciences, Shiraz University, Shiraz 71454, Iran
E-mail: m.mehdiloozad@gmail.com

**Kaoru Tone**
National Graduate Institute for Policy Studies, Tokyo, Japan

**Mohammad Bagher Ahmadi**
Department of Mathematics, College of Sciences, Shiraz University, Shiraz 71454, Iran

# ABSTRACT

The *strict complementary slackness condition (SCSC)* is an important concept in the duality theory of linear programming (LP). The current study aims at extending this concept to the framework of linear fractional programming (LFP). First, we define this concept in this framework and demonstrate the existence of a *strict complementary solution* — a pair of primal and dual optimal solutions satisfying the SCSC. Second, we show that the problem of finding such a solution reduces to that of identifying a relative interior point of a polyhedron. More recently, Mehdiloozad et al. (2016) have addressed the latter problem by proposing an LP problem. Using their proposed LP problem, we finally develop two procedures for finding a strict complementary solution.

# 1. INTRODUCTION

A mathematical programming problem is specified as a (primal) *linear fractional programming (LFP)* when a linear fractional objective is optimized subject to a set of linear constraints. The LFP problem has attracted considerable research interests due to its usefulness in real-world applications. Details on the applications of the LFP problem can be found in Bajalinov (2003). In accordance with the duality theory of mathematical programming, a problem, called the *dual LFP problem*, can be associated with the primal LFP problem. In the literature, different types of the dual LFP problem have been proposed (Bajalinov 2003), one of which has the form of a linear programming (LP) problem. Chadha (1971) formulated the LP form of the dual LFP problem and proved some duality results *directly*. In his formulation, constant terms appear in neither of the numerator and denominator of the objective function of the primal LFP problem. Subsequently, Chadha and Chadha (2007) extended his work to the general case where the constant terms are taken into account. As a main result, they showed that the *complementary slackness condition (CSC)* holds true between the primal and dual LFP problems.

As is well known, an effective method for solving the LFP problem is to transform it into an equivalent LP problem by using the Charnes–Cooper (1962) transformation. With a little care, it can be observed that the dual LFP problem of Chadha and Chadha (2007) is nothing else than the dual of the transformed LP problem. This indicates that all the duality relationships between the primal and dual LFP problems can be *indirectly* deduced from the duality of linear programming. Based on this fact, we aim at extending the Goldman–Tucker (1956) theorem to the case of LFP problems.

The Goldman–Tucker theorem is an important duality result in linear programming that establishes the *strict complementary slackness condition (SCSC)*. The SCSC asserts the existence of an optimal primal–dual solution for which the sum of each pair of complementary variables is positive. A feasible primal–dual solution satisfying the SCSC, for brevity, is called a *strict complementary solution*. For more details about the SCSC, the interested readers can refer to the recent work by Mehdiloozad et al. (2016).

To the best of our knowledge, no research study has been conducted to date on examining the concept of SCSC in the framework of linear fractional programming. Motivated by this, we formally define this concept in this framework and demonstrate the existence of a strict complementary solution by means of the Goldman–Tucker theorem. Under the existence of a strict complementary solution, we then deal with the problem of finding such a solution or, equivalently, the problem of identifying a relative interior point of a polyhedron. Mehdiloozad et al. (2016) have addressed the latter problem by proposing an LP problem. Using their proposed LP problem, we develop two procedures for finding a strict complementary solution.

The remainder of this paper unfolds as follows. Section 2 provides some basic concepts related to the duality in linear fractional programming. Section 3 introduces and discusses the notion of SCSC. It then develops two methods for finding a strict complementary solution. Section 4 illustrates our proposed approach with a numerical example. Section 5 concludes with some remarks.

## 2. Duality in linear fractional programming

As far as notations are concerned, the *n*-dimensional Euclidean space is denoted by $\mathbb{R}^n$ with non-negative orthant denoted by $\mathbb{R}^n_+$. We symbolize the sets by capital letters and their members by lower-case letters. We also denote vectors and matrices in bold letters, vectors in lower case and matrices in upper case. All vectors are column vectors. We denote the transpose of vectors and matrices by a superscript $T$. We also use $\mathbf{0}_n$ and $\mathbf{1}_n$ to show *n*-dimensional vectors with the values of 0 and 1 in every entry, respectively. Furthermore, we denote the $n \times n$ identity matrix by $\mathbf{I}_n$. For any $\mathbf{x}, \mathbf{y} \in \mathbb{R}^n$, we write $\mathbf{x} \leq \mathbf{y}$ to indicate that $x_j \leq y_j$ for all $j = 1,...,n$.

A *linear fractional programming* (LFP) program is a mathematical programming problem that optimizes a linear fractional objective subject to a set of linear constraints. Formally, the primal LFP problem is of the following form:

$$\max \quad f(\mathbf{x}) = \frac{\mathbf{c}^T \mathbf{x} + \alpha}{\mathbf{d}^T \mathbf{x} + \beta}$$
$$\text{subject to} \qquad (1)$$
$$\mathbf{x} \in X = \{\mathbf{x} \mid \mathbf{A}\mathbf{x} \leq \mathbf{b}, \ \mathbf{x} \geq \mathbf{0}_n\},$$

where $\mathbf{A}$ is an $m \times n$ matrix, $\mathbf{c}$, $\mathbf{d}$ are $n \times 1$ vectors, $\mathbf{b}$ is an $m \times 1$ vector, and $\alpha$, $\beta$ are scalars. It is assumed that $X$ is regular (non-empty and bounded), the objective function $f$ is continuously differentiable, and $\mathbf{d}^T \mathbf{x} + \beta > 0$ for all $\mathbf{x} \in X$.

Using the Charnes–Cooper (1962) transformation, we transform problem (1) into an equivalent LP problem. For this, let us define

$$t := \frac{1}{\mathbf{d}^T \mathbf{x} + \beta} \quad \& \quad \bar{\mathbf{x}} := t\mathbf{x}. \qquad (2)$$

Then, by multiplying both sides of the constraints $\mathbf{A}\mathbf{x} \leq \mathbf{b}$ by $t$, problem (1) is rewritten equivalently as the following LP problem:

$$\begin{aligned}&\max\quad \mathbf{c}^T\overline{\mathbf{x}}+\alpha t\\&\text{subject to}\\&\begin{pmatrix}\overline{\mathbf{x}}\\t\end{pmatrix}\in \overline{X}=\left\{\begin{pmatrix}\overline{\mathbf{x}}\\t\end{pmatrix}\bigg|\ \mathbf{A}\overline{\mathbf{x}}\le \mathbf{b}t,\ \mathbf{d}^T\overline{\mathbf{x}}+\beta t=1,\ \overline{\mathbf{x}}\ge \mathbf{0}_n,\ t\ge 0\right\}.\end{aligned} \quad (3)$$

Since $t>0$ for all $\begin{pmatrix}\overline{\mathbf{x}}\\t\end{pmatrix}\in \overline{X}$, the Charnes–Cooper transformation establishes a one-to-one correspondence between problems (1) and (3). More precisely, $\frac{1}{t}\overline{\mathbf{x}}\in X$ if and only if $\begin{pmatrix}\overline{\mathbf{x}}\\t\end{pmatrix}\in \overline{X}$, and the statement holds also true at optimality.

Now, we consider the dual of problem (3) given by

$$\begin{aligned}&\min\quad g(\mathbf{y},z)=z\\&\text{subject to}\\&\begin{pmatrix}\mathbf{y}\\z\end{pmatrix}\in Y=\left\{\begin{pmatrix}\mathbf{y}\\z\end{pmatrix}\bigg|\ \mathbf{A}^T\mathbf{y}+\mathbf{d}z\ge \mathbf{c},\ -\mathbf{b}^T\mathbf{y}+\beta z=\alpha,\ \mathbf{y}\ge \mathbf{0}_m\right\}.\end{aligned} \quad (4)$$

**Remark 2.1** Note that the inequality constraint $-\mathbf{b}^T\mathbf{y}+\beta z\ge \alpha$ is written in equality form in problem (4) since its corresponding primal variable is positive at optimality, i.e., $t^*>0$.

In accordance with Chadha and Chadha (2007), we define problem (4) as the *dual LFP problem* and call, for brevity, a pair of primal and dual feasible and optimal solutions a *feasible primal-dual solution* and an *optimal primal-dual solution*, respectively. Then, the following four theorems demonstrate the duality relationships between problems (1) and (4).

**Theorem 2.1** (*Weak duality*) $f(\mathbf{x})\le g(\mathbf{y},z)$ holds for any feasible primal-dual solution $(\mathbf{x},(\mathbf{y},z))$.

**Theorem 2.2** (*Optimality criteria*) Let $\left(\mathbf{x}^*,\left(\mathbf{y}^*,z^*\right)\right)$ is a feasible primal-dual solution such that $f\left(\mathbf{x}^*\right)=g\left(\mathbf{y}^*,z^*\right)$. Then, it is optimal.

**Theorem 2.3** (*Strong duality*) If $\mathbf{x}^*$ optimizes problem (1), then there exists $\left(\mathbf{y}^*,z^*\right)$ that optimizes problem (4) and $f\left(\mathbf{x}^*\right)=g\left(\mathbf{y}^*,z^*\right)$.

**Theorem 2.4** (*complementary slackness*) $\left(\mathbf{x}^*,\left(\mathbf{y}^*,z^*\right)\right)$ is an optimal primal-dual solution if and only if

$$\mathbf{x}^{*T}\mathbf{v}^* = 0 \quad \& \quad \mathbf{y}^{*T}\mathbf{u}^* = 0, \tag{5}$$

where $\mathbf{u}^* := \mathbf{b} - \mathbf{A}\mathbf{x}^*$ and $\mathbf{v}^* := \mathbf{A}^T\mathbf{y}^* + \mathbf{d}z^* - \mathbf{c}$ are the optimal primal and dual slack vectors, respectively.

Equation (5) is known as the *complementary slackness condition (CSC)*, and $\left(x_j^*, v_j^*\right)$, $j=1,...,n$, and $\left(y_i^*, u_i^*\right)$, $i=1,...,m$, are called as the *complementary variables*.

## 3. THE SCSC IN LINEAR FRACTIONAL PROGRAMMING

### 3.1. DEFINITION OF THE SCSC

While the CSC asserts that the multiplication of each pair of complementary variables is zero, it does not guarantee the positivity of the sum of these variables. For a primal–dual LFP optimal solutions, if both complementarity and positivity conditions hold, then exactly one of the complementary variables takes a positive value and the other one is zero. Consequently, the following unique partitions, called *optimal partitions*, can be induced for the index sets $\{1,...,n\}$ and $\{1,...,m\}$:

$$\{j \mid x_j^* > 0\} \cup \{j \mid v_j^* > 0\} = \{1,...,n\}, \quad \{j \mid x_j^* > 0\} \cap \{j \mid v_j^* > 0\} = \emptyset,$$
$$\{i \mid u_i^* > 0\} \cup \{i \mid y_i^* > 0\} = \{1,...,m\}, \quad \{i \mid u_i^* > 0\} \cap \{i \mid y_i^* > 0\} = \emptyset. \quad (6)$$

Now, a question arises as to how determine these partitions. To address this question, first, we formally express the satisfactions of both complementarity and positivity conditions in the following definition.

**Definition 3.1** For a feasible primal–dual solution $\left((\mathbf{x}^*, \mathbf{u}^*), (\mathbf{y}^*, z^*, \mathbf{v}^*)\right)$, the *strict complementary slackness condition (SCSC)* holds if and only if it satisfies the following condition as well as condition (5):

$$\mathbf{x}^* + \mathbf{v}^* > \mathbf{0}_n \quad \& \quad \mathbf{y}^* + \mathbf{u}^* > \mathbf{0}_m. \quad (7)$$

Then, the given pair is called a *strict complementary solution*.

**Theorem 3.1** There exists a strict complementary solution.

*Proof* From the Goldman–Tucker (1956) theorem, there exists a strict complementary solution for the LP problems (3) and (4), namely $\left((\bar{\mathbf{x}}^*, t^*, \bar{\mathbf{u}}^*), (\mathbf{y}^*, z^*, \mathbf{v}^*)\right)$, where $\bar{\mathbf{u}}^* := \mathbf{b}t^* - \mathbf{A}\bar{\mathbf{x}}^*$ and $\mathbf{v}^* := \mathbf{A}^T \mathbf{y}^* + \mathbf{d}z^* - \mathbf{c}$. Then, we have

$$\bar{\mathbf{x}}^{*T}\mathbf{v}^* = \mathbf{y}^{*T}\bar{\mathbf{u}}^* = 0 \quad \& \quad \bar{\mathbf{x}}^* + \mathbf{v}^* > \mathbf{0}_n \quad \& \quad \mathbf{y}^* + \bar{\mathbf{u}}^* > \mathbf{0}_m. \quad (8)$$

Since $t^* > 0$, it is immediate by (2) that $(\mathbf{x}^*, \mathbf{u}^*) := \frac{1}{t^*}(\bar{\mathbf{x}}^*, \bar{\mathbf{u}}^*)$ and $(\mathbf{y}^*, z^*, \mathbf{v}^*)$ are optimal solutions to problems (1) and (4), respectively, for which conditions (5) and (7) hold. This completes the proof. ∎

Theorem 3.1 demonstrates the existence of a strict complementary solution. However, it does not specify how to identify such a solution. Hence, we proceed to develop a method for dealing with this issue.

To find a strictly complementary solution, we include, respectively, the slack vectors $\mathbf{u} \in \mathbb{R}_+^m$ and $\mathbf{v} \in \mathbb{R}_+^n$ into the inequality constraints of problems (3) and (4) as

$$\max \quad \mathbf{c}^T \bar{\mathbf{x}} + \alpha t$$
$$\text{subject to} \tag{9}$$
$$\begin{pmatrix} \bar{\mathbf{x}} \\ t \\ \mathbf{u} \end{pmatrix} \in S_P = \left\{ \begin{pmatrix} \bar{\mathbf{x}} \\ t \\ \mathbf{u} \end{pmatrix} \middle| \mathbf{A}\bar{\mathbf{x}} + \mathbf{u} = \mathbf{b}t,\ \mathbf{d}^T\bar{\mathbf{x}} + \beta t = 1,\ \bar{\mathbf{x}} \geq \mathbf{0}_n,\ t \geq 0,\ \mathbf{u} \geq \mathbf{0}_m \right\}.$$

$$\min \quad g(\mathbf{y}, z) = z$$
$$\text{subject to} \tag{10}$$
$$\begin{pmatrix} \mathbf{y} \\ z \\ \mathbf{v} \end{pmatrix} \in S_D = \left\{ \begin{pmatrix} \mathbf{y} \\ z \\ \mathbf{v} \end{pmatrix} \middle| \mathbf{A}^T\mathbf{y} + \mathbf{d}z - \mathbf{v} = \mathbf{c},\ -\mathbf{b}^T\mathbf{y} + \beta z = \alpha,\ \mathbf{y} \geq \mathbf{0}_m,\ \mathbf{v} \geq \mathbf{0}_n \right\}.$$

Let $F_P^*$ and $F_D^*$ denote the optimal faces of problems (9) and (10), respectively. If the optimal objective value of problems (9) and (10) is equal to $\theta^*$, then $F_P^*$ and $F_D^*$ can be formulated as follows:

$$F_P^* = \left\{ \begin{pmatrix} \bar{\mathbf{x}} \\ t \\ \bar{\mathbf{u}} \end{pmatrix} \middle| \begin{bmatrix} \mathbf{A} & -\mathbf{b} & \mathbf{I}_m \\ \mathbf{d}^T & \beta & \mathbf{0}_m^T \\ \mathbf{c}^T & \alpha & \mathbf{0}_m^T \end{bmatrix} \begin{pmatrix} \bar{\mathbf{x}} \\ t \\ \bar{\mathbf{u}} \end{pmatrix} = \begin{pmatrix} \mathbf{0}_m \\ 1 \\ \theta^* \end{pmatrix},\ \begin{pmatrix} \bar{\mathbf{x}} \\ t \\ \bar{\mathbf{u}} \end{pmatrix} \geq \mathbf{0}_{m+n+1} \right\}, \tag{11}$$

$$F_D^* = \left\{ \begin{pmatrix} \mathbf{y} \\ z \\ \mathbf{v} \end{pmatrix} \middle| \begin{bmatrix} \mathbf{A}^T & \mathbf{d} & -\mathbf{I}_n \\ -\mathbf{b}^T & \beta & \mathbf{0}_n^T \\ \mathbf{0}_m^T & 1 & \mathbf{0}_n^T \end{bmatrix} \begin{pmatrix} \mathbf{y} \\ z \\ \mathbf{v} \end{pmatrix} = \begin{pmatrix} \mathbf{c} \\ \alpha \\ \theta^* \end{pmatrix},\ \begin{pmatrix} \mathbf{y} \\ z \\ \mathbf{v} \end{pmatrix} \geq \mathbf{0}_{m+n+1} \right\}. \tag{12}$$

We now demonstrate how the SCSC relates to the relative interiors of these two optimal faces, denoted by $ri(F_P^*)$ and $ri(F_D^*)$.

**Theorem 3.1** Let $\begin{pmatrix} \bar{\mathbf{x}}^* \\ t^* \\ \bar{\mathbf{u}}^* \end{pmatrix} \in ri(F_P^*)$ and $\begin{pmatrix} \mathbf{y}^* \\ z^* \\ \mathbf{v}^* \end{pmatrix} \in ri(F_D^*)$. Then, $\left( \frac{1}{t^*}(\bar{\mathbf{x}}^*, \bar{\mathbf{u}}^*), (\mathbf{y}^*, z^*, \mathbf{v}^*) \right)$ is a strict complementary solution.

*Proof* The proof readily follows from Theorem 5.1 in Mehdiloozad et al. (2016) by using the transformation (2) and hence is omitted. ∎

### 3.2 FINDING A STRICT COMPLEMENTARY SOLUTION

Let us consider the polyhedron $P$ defined by

$$P = \{\mathbf{x} | \mathbf{A}\mathbf{x} = \mathbf{b},\ \mathbf{x} \geq \mathbf{0}_n\}, \tag{13}$$

where $\mathbf{A}$ is a matrix of order $m \times n$ and $\mathbf{b} \in \mathbb{R}^m$. Then, according to Theorem 3.1, the problem of finding a strict complementary solution reduces to the problem of identifying a relative interior point of $P$. On the other hand, Theorem 4.1 in Mehdiloozad et al. (2016) follows that this problem is equivalent to the problem of finding a maximal element of $P$ ($\mathbf{x}^{max}$) — an element with the maximum number of positive components. To address the latter problem, Mehdiloozad et al. (2016) have developed the following LP problem:

$$\begin{aligned} \max \quad & \mathbf{1}_n^T \mathbf{x}^2 + w^2 \\ \text{subject to} \quad & \\ & [\mathbf{A}, -\mathbf{b}] \begin{bmatrix} \mathbf{x}^1 + \mathbf{x}^2 \\ w^1 + w^2 \end{bmatrix} = \mathbf{0}_m, \\ & \mathbf{0}_{n+1} \leq \begin{pmatrix} \mathbf{x}^1 \\ w^1 \end{pmatrix},\ \mathbf{0}_{n+1} \leq \begin{pmatrix} \mathbf{x}^2 \\ w^2 \end{pmatrix} \leq \mathbf{1}_{n+1}. \end{aligned} \tag{14}$$

Moreover, they have demonstrated that $\mathbf{x}^{max}$ can be obtained as follows:

$$\mathbf{x}^{max} = \frac{1}{w^{1*} + w^{2*}} \left( \mathbf{x}^{1*} + \mathbf{x}^{2*} \right). \tag{15}$$

### 3.2.1 FIRST APPROACH

In this approach, we assume that the optimal objective value of problems (9) and (10) is known and is equal to $\theta^*$. Then, by Theorem 3.1 and problem (14), we propose the following LP problems to identify relative interior points of $F_P^*$ and $F_D^*$:

$$\max \quad \sum_{j=1}^{n} \bar{x}_j^2 + \sum_{i=1}^{m} \bar{u}_i^2 + w_P^2$$

subject to

$$\begin{bmatrix} \mathbf{A} & -\mathbf{b} & \mathbf{I}_m & \mathbf{0}_m \\ \mathbf{d}^T & \beta & \mathbf{0}_m^T & -1 \\ \mathbf{c}^T & \alpha & \mathbf{0}_m^T & -\theta^* \end{bmatrix} \begin{bmatrix} \bar{\mathbf{x}}^1 + \bar{\mathbf{x}}^2 \\ p \\ \bar{\mathbf{u}}^1 + \bar{\mathbf{u}}^2 \\ w_P^1 + w_P^2 \end{bmatrix} = \mathbf{0}_{m+2}, \qquad (16)$$

$$p \geq 0, \quad \begin{pmatrix} \bar{\mathbf{x}}^1 \\ \bar{\mathbf{u}}^1 \\ w_P^1 \end{pmatrix} \geq \mathbf{0}_{m+n+1}, \quad \mathbf{0}_{m+n+1} \leq \begin{pmatrix} \bar{\mathbf{x}}^2 \\ \bar{\mathbf{u}}^2 \\ w_P^2 \end{pmatrix} \leq \mathbf{1}_{m+n+1}.$$

$$\max \quad \sum_{i=1}^{m} y_i^2 + \sum_{j=1}^{n} v_j^2 + w_D^2$$

subject to

$$\begin{bmatrix} \mathbf{A}^T & \mathbf{d} & -\mathbf{I}_n & -\mathbf{c} \\ -\mathbf{b}^T & \beta & \mathbf{0}_n^T & -\alpha \\ \mathbf{0}_m^T & 1 & \mathbf{0}_n^T & -\theta^* \end{bmatrix} \begin{bmatrix} \mathbf{y}^1 + \mathbf{y}^2 \\ q \\ \mathbf{v}^1 + \mathbf{v}^2 \\ w_D^1 + w_D^2 \end{bmatrix} = \mathbf{0}_{n+2}, \qquad (17)$$

$$\begin{pmatrix} \mathbf{y}^1 \\ \mathbf{v}^1 \\ w_D^1 \end{pmatrix} \geq \mathbf{0}_{m+n+1}, \quad \mathbf{0}_{m+n+1} \leq \begin{pmatrix} \mathbf{y}^2 \\ \mathbf{v}^2 \\ w_D^2 \end{pmatrix} \leq \mathbf{1}_{m+n+1}.$$

Now, let $\left(p^*, \bar{\mathbf{x}}^{1*}, \bar{\mathbf{x}}^{2*}, \bar{\mathbf{u}}^{1*}, \bar{\mathbf{u}}^{2*}, w_P^{1*}, w_P^{2*}\right)$ and $\left(q^*, \mathbf{y}^{1*}, \mathbf{y}^{2*}, \mathbf{v}^{1*}, \mathbf{v}^{2*}, w_D^{1*}, w_D^{2*}\right)$ be optimal solutions to problems (16) and (17), respectively. Then,

$$\begin{pmatrix} \bar{\mathbf{x}}^* \\ t^* \\ \bar{\mathbf{u}}^* \end{pmatrix} := \frac{1}{w_P^{1*} + w_P^{2*}} \begin{pmatrix} \bar{\mathbf{x}}^{1*} + \bar{\mathbf{x}}^{2*} \\ p^* \\ \bar{\mathbf{u}}^{1*} + \bar{\mathbf{u}}^{2*} \end{pmatrix} \in ri\left(F_P^*\right), \qquad (18)$$

$$\begin{pmatrix} \mathbf{y}^* \\ z^* \\ \mathbf{v}^* \end{pmatrix} := \frac{1}{w_D^{1*} + w_D^{2*}} \begin{pmatrix} \mathbf{y}^{1*} + \mathbf{y}^{2*} \\ q^* \\ \mathbf{v}^{1*} + \mathbf{v}^{2*} \end{pmatrix} \in ri\left(F_D^*\right). \tag{19}$$

Therefore, as per Theorem 3.1, a strictly complementary solution can be obtained from (18) and (19).

### 3.2.2 SECOND APPROACH

Our first approach of identifying a strict complementary solution requires the knowledge of the optimal objective value of problems (9) and (10). Now, we aim at proposing an alternative approach that is exempt from this requirement. Toward this, we exploit the fact that the optimality of a feasible primal-dual solution in linear programming follows from the equality of the primal and dual objective function values. Based on this fact, we define

$$F_{PD}^* := \left\{ \begin{pmatrix} \overline{\mathbf{x}} \\ t \\ \overline{\mathbf{u}} \\ \mathbf{y} \\ z \\ \mathbf{v} \end{pmatrix} \middle| \begin{bmatrix} \mathbf{A} & -\mathbf{b} & \mathbf{I}_m & \mathbf{0}_{m \times m} & 0 & \mathbf{0}_{m \times n} \\ \mathbf{d}^T & \beta & \mathbf{0}_m^T & \mathbf{0}_m^T & 0 & \mathbf{0}_n^T \\ \mathbf{0}_{n \times n} & 0 & \mathbf{0}_{n \times m} & \mathbf{A}^T & \mathbf{d} & -\mathbf{I}_n \\ \mathbf{0}_n^T & 0 & \mathbf{0}_m^T & -\mathbf{b}^T & \beta & \mathbf{0}_n^T \\ \mathbf{c}^T & \alpha & \mathbf{0}_m^T & \mathbf{0}_m^T & -1 & \mathbf{0}_n^T \end{bmatrix} \begin{bmatrix} \overline{\mathbf{x}} \\ t \\ \overline{\mathbf{u}} \\ \mathbf{y} \\ z \\ \mathbf{v} \end{bmatrix} = \begin{bmatrix} \mathbf{0}_m \\ 1 \\ \mathbf{c} \\ \alpha \\ 0 \end{bmatrix}, \begin{pmatrix} \overline{\mathbf{x}} \\ t \\ \overline{\mathbf{u}} \\ \mathbf{y} \\ z \\ \mathbf{v} \end{pmatrix} \geq \mathbf{0}_{2m+2n+2} \right\}. \tag{20}$$

The set of constraints defining $F_{PD}^*$ consists of the sets of constraints of problems (9) and (10) together with $\mathbf{c}^T \overline{\mathbf{x}} + \alpha t - z = 0$. Therefore, for any $\begin{pmatrix} \overline{\mathbf{x}} \\ t \\ \overline{\mathbf{u}} \\ \mathbf{y} \\ z \\ \mathbf{v} \end{pmatrix} \in F_{PD}^*$, we have $\begin{pmatrix} \overline{\mathbf{x}} \\ t \\ \overline{\mathbf{u}} \end{pmatrix} \in F_P^*$ and $\begin{pmatrix} \mathbf{y} \\ z \\ \mathbf{v} \end{pmatrix} \in F_D^*$.

Hence, a maximal element of $F_{PD}^*$ determines a strict complementary solution. In view of this fact, we formulate the following LP problem for finding a strict complementary solution:

$$\max \quad \sum_{j=1}^{n} x_j^2 + \sum_{i=1}^{m} u_i^2 + \sum_{i=1}^{m} y_i^2 + \sum_{j=1}^{n} v_j^2 + w^2$$

subject to

$$\begin{bmatrix} \mathbf{A} & -\mathbf{b} & \mathbf{I}_m & \mathbf{0}_{m \times m} & 0 & \mathbf{0}_{m \times n} & \mathbf{0}_m \\ \mathbf{d}^T & \beta & \mathbf{0}_m^T & \mathbf{0}_m^T & 0 & \mathbf{0}_n^T & 1 \\ \mathbf{0}_{n \times n} & 0 & \mathbf{0}_{n \times m} & \mathbf{A}^T & \mathbf{d} & -\mathbf{I}_n & \mathbf{c} \\ \mathbf{0}_n^T & 0 & \mathbf{0}_m^T & -\mathbf{b}^T & \beta & \mathbf{0}_n^T & \alpha \\ \mathbf{c}^T & \alpha & \mathbf{0}_m^T & \mathbf{0}_m^T & -1 & \mathbf{0}_n^T & 0 \end{bmatrix} \begin{bmatrix} \overline{\mathbf{x}}^1 + \overline{\mathbf{x}}^2 \\ p \\ \overline{\mathbf{u}}^1 + \overline{\mathbf{u}}^2 \\ \mathbf{y}^1 + \mathbf{y}^2 \\ q \\ \mathbf{v}^1 + \mathbf{v}^2 \\ w^1 + w^2 \end{bmatrix} = \mathbf{0}_{2m+2n+2},$$

$$p \geq 0, \begin{pmatrix} \mathbf{x}^1 \\ \overline{\mathbf{u}}^1 \\ \mathbf{y}^1 \\ \mathbf{v}^1 \\ w^1 \end{pmatrix} \geq \mathbf{0}_{2m+2n+1}, \quad \mathbf{0}_{2m+2n+1} \leq \begin{pmatrix} \mathbf{x}^2 \\ \overline{\mathbf{u}}^2 \\ \mathbf{y}^2 \\ \mathbf{v}^2 \\ w^2 \end{pmatrix} \leq \mathbf{1}_{2m+2n+1}. \tag{21}$$

Let $\left( p^*, q^*, \overline{\mathbf{x}}^{1*}, \overline{\mathbf{x}}^{2*}, \overline{\mathbf{u}}^{1*}, \overline{\mathbf{u}}^{2*}, \mathbf{y}^{1*}, \mathbf{y}^{2*}, \mathbf{v}^{1*}, \mathbf{v}^{2*}, w^{1*}, w^{2*} \right)$ be an optimal solution to problem (21). Then,

$$\begin{pmatrix} \overline{\mathbf{x}}^* \\ t^* \\ \overline{\mathbf{u}}^* \end{pmatrix} := \frac{1}{w^{1*}+1} \begin{pmatrix} \overline{\mathbf{x}}^{1*} + \overline{\mathbf{x}}^{2*} \\ p^* \\ \overline{\mathbf{u}}^{1*} + \overline{\mathbf{u}}^{2*} \end{pmatrix} \in ri\left( F_P^* \right), \tag{22}$$

$$\begin{pmatrix} \mathbf{y}^* \\ z^* \\ \mathbf{v}^* \end{pmatrix} := \frac{1}{w^{1*}+1} \begin{pmatrix} \mathbf{y}^{1*} + \mathbf{y}^{2*} \\ q^* \\ \mathbf{v}^{1*} + \mathbf{v}^{2*} \end{pmatrix} \in ri\left( F_D^* \right). \tag{23}$$

Then, a strictly complementary solution is obtained using Theorem 3.1.

## 4. NUMERICAL EXAMPLE

Let us consider the following primal LFP problem

$$\max \quad \frac{6x_1 + 3x_2 + 6}{5x_1 + 2x_2 + 5}$$

subject to

$$\begin{bmatrix} 2 & 1 \\ -2 & 1 \end{bmatrix} \begin{bmatrix} x_1 \\ x_2 \end{bmatrix} + \begin{bmatrix} u_1 \\ u_2 \end{bmatrix} = \begin{bmatrix} 6 \\ 2 \end{bmatrix}, \quad (24)$$

$$\begin{pmatrix} x_1 \\ x_2 \end{pmatrix}, \begin{pmatrix} u_1 \\ u_2 \end{pmatrix} \geq \begin{pmatrix} 0 \\ 0 \end{pmatrix},$$

whose associated dual LFP problem is

$$\min \quad z$$

subject to

$$\begin{bmatrix} 2 & -2 \\ 1 & 1 \end{bmatrix} \begin{bmatrix} y_1 \\ y_2 \end{bmatrix} + \begin{bmatrix} 5 \\ 2 \end{bmatrix} z - \begin{bmatrix} v_1 \\ v_2 \end{bmatrix} = \begin{bmatrix} 6 \\ 3 \end{bmatrix}, \quad (25)$$

$$-\begin{bmatrix} 6 & 2 \end{bmatrix} \begin{bmatrix} y_1 \\ y_2 \end{bmatrix} + 5z = 6,$$

$$\begin{pmatrix} x_1 \\ x_2 \end{pmatrix}, \begin{pmatrix} u_1 \\ u_2 \end{pmatrix} \geq \begin{pmatrix} 0 \\ 0 \end{pmatrix}.$$

Then, by the graphical method described in Bajalinov (2003), it can be easily verified that the optimal objective value is equal to 1.333 and all the primal optimal solutions are the convex combinations of $\begin{pmatrix} 0 \\ 2 \\ 4 \\ 0 \end{pmatrix}$ and $\begin{pmatrix} 1 \\ 4 \\ 0 \\ 0 \end{pmatrix}$. Moreover, using the duality results, all the optimal primal-dual solutions are determined as

$$\begin{pmatrix} x_1^* \\ x_2^* \\ u_1^* \\ u_2^* \end{pmatrix} = \begin{pmatrix} 1-\alpha \\ 4-2\alpha \\ 4\alpha \\ 0 \end{pmatrix}, \quad \alpha \in [0,1] \quad \& \quad \begin{pmatrix} v_1^* \\ v_2^* \\ y_1^* \\ y_2^* \end{pmatrix} = \begin{pmatrix} 0 \\ 0 \\ 0 \\ 0.333 \end{pmatrix}. \quad (26)$$

It is obvious that $\left(\left(\mathbf{x}^*, \mathbf{u}^*\right), \left(\mathbf{y}^*, z^*, \mathbf{v}^*\right)\right)$ is a strict complementary solution for $\alpha \in (0,1)$. Then, $\sigma(\mathbf{x}^*) = \{1,2\}$ and $\sigma(\mathbf{v}^*) = \emptyset$ constitute the optimal partition of the index set $\{1,2\}$ ($n = 2$). Moreover, $\sigma(\mathbf{u}^*) = \{1\}$ and $\sigma(\mathbf{y}^*) = \{2\}$ build the optimal partition of the index set $\{1,2\}$ ($m = 2$). If $\alpha = 0$ or $\alpha = 1$, however, the optimal primal-dual solution is not strict complementary. This is because in either of these two cases, the sum of each pair of the complementary variables is not positive.

Applying each of our proposed approaches, we find the following strict complementary solution:

$$\begin{pmatrix} x_1^* \\ x_2^* \\ u_1^* \\ u_2^* \end{pmatrix} = \begin{pmatrix} 0.8 \\ 3.6 \\ 0.8 \\ 0 \end{pmatrix} \quad \& \quad \begin{pmatrix} v_1^* \\ v_2^* \\ y_1^* \\ y_2^* \end{pmatrix} = \begin{pmatrix} 0 \\ 0 \\ 0 \\ 0.333 \end{pmatrix}, \tag{27}$$

which is corresponding to $\alpha = 0.2$. We have carried out all the computations by developing a computer program, given in Appendix A, using the GAMS optimization software.

## 5. CONCLUSIONS

As an effective technique, the Charnes-Cooper transformation solves the primal LFP problem by converting it into an equivalent LP problem. By defining the dual of the transformed LP problem as the dual LFP problem, the duality of linear programming can be used for demonstrating the duality results between the primal and dual LFP problems. Based on this fact, we introduced the concept of SCSC in linear fractional programming. Then, by virtue of the Goldman-Tucker theorem, we proved the existence of a strict complementary solution. Then, we developed two approaches with different strategies for finding a strict complementary solution.

Our first and second approaches identify a strictly complementary solution by solving two and one LP problem(s), respectively. Therefore, they make straightforward to apply the ordinary simplex algorithm of linear programming for identifying a strictly complementary solution. As regards the preference on using the proposed approaches, note that each of the LP problems of our first approach has less number of constraints than the LP problem of our second approach. In view of this fact, the first one is particularly recommended in situations where only primal or dual part of a strictly complementary solution needs to be found. For example, in the theory of data envelopment analysis (Charnes et al. 1987), the global reference set[1] of an inefficient decision making unit can be determined using only the primal part of strict complementary solution to the SBM model of Tone (2001).

In linear programming, an interesting method to find a strict complementary solution is to apply the so-called Balinski–Tucker (optimal) tableau (Balinski and Tucker 1969). Hence, using the Charnes–Cooper transformation, this method is able to indirectly generate a strict complementary solution from the transformed LP form of the primal LFP and the dual LFP problem. Nonetheless, as a future research subject, we suggest modifying this method for *direct* identification of a strict complementary solution in linear fractional programming.

# APPENDIX A

The computer program (written in GAMS) for identifying a strictly complementary solution:

```
1   Sets
2       i       row number of matrix A    /i1*i2/
3       j       column number of matrix A /j1*j2/;
4
5   Table A(i,j)
6           j1      j2
7   i1      2       1
8   i2      -2      1 ;
9
10  Parameters
11      b/i1    6
12       i2     2   /
13      c/j1    6
14       j2     3   /
15      d/j1    5
16       j2     2   / ;
17
```

---

[1] For more details about the concept of global reference set, the interested readers are referred to Mehdiloozad et al. (2016) and Mehdiloozad (2016).

```
18  Scalars
19      Alpha
20      Beta ;
21
22  Alpha=6;
23  Beta =5;
24
25  File Prog / Result.txt/;
26  Put Prog;
27
28  *************************************************************
29  *Stage 1: Solving problem (3)
30  ****************************
31
32  Free Variables
33      Theta;
34
35  Positive Variables
36      xbar(j)
37      t ;
38
39  Scalar
40      ThetaStar ;
41
42  Equations
43      Obj
44        Con1
45        Con2 ;
46
47      Obj..          Theta =E= Sum(j, c(j)*xbar(j)) + Alpha*t;
48      Con1(i)..      Sum(j, a(i,j)*xbar(j))         =L= b(i)*t;
49      Con2..         Sum(j,  d(j)*xbar(j)) + Beta*t =E= 1;
50
51  Model MainLFP     / Obj, Con1, Con2 /;
52
53  Put /'Finding the Optimal Obj. Value (ThetaStar)';
54  Put /'-----------------------------------------'/;
55  Option LP=CONOPT;
56  Solve MainLFP using LP Maximizing Theta;
57      Put 'Obj = ':>6; Put Theta.L:<10:3;
58      ThetaStar=Theta.L;
59  Put /'-----------------------------------------'/;
60
61  ****************************
62  *Stage 1: Solving problem (3)
63  *************************************************************
64
```

```
65   *************************************************************
66   *First approach
67   ***************
68
69   Positive Variables
70        xbar1(j)
71        xbar2(j)
72        ubar1(i)
73        ubar2(i)
74        y1(i)
75        y2(i)
76        v1(j)
77        v2(j)
78        w1
79        w2
80        p ;
81
82   Free variable
83        q ;
84
85   xbar2.up(j) = 1;
86   ubar2.up(i) = 1;
87   y2.up(i) = 1;
88   v2.up(j) = 1;
89   w2.up = 1;
90
91   Parameters
92        XbarStar
93        tStar
94        UbarStar
95        YStar(i)
96        zStar
97        VStar(j) ;
98
99   Equations
100       ObjP
101         ConP1
102         ConP2
103         ConP3
104       ObjD
105         ConD1
106         ConD2
107         ConD3   ;
108
109       ObjP..         Theta =E= Sum(j, xbar2(j)) + Sum(i, ubar2(i)) + w2;
110       ConP1(i)..       Sum(j, a(i,j)*(xbar1(j)+xbar2(j))) - b(i)*p + ubar1(i)+ubar2(i)
=E= 0;
```

```
111         ConP2..             Sum(j, d(j)*(xbar1(j)+xbar2(j))) + Beta*p - w1+w2 =E= 0;
112         ConP3..             Sum(j, c(j)*(xbar1(j)+xbar2(j))) + Alpha*p - (w1+w2)*ThetaStar
=E= 0;
113
114       ObjD..        Theta =E= Sum(i, y2(i)) + Sum(j, v2(j)) + w2;
115         ConD1(j)..      Sum(i, a(i,j)*(y1(i)+y2(i))) + d(j)*q - v1(j)-v2(j) - c(j)*
(w1+w2) =E= 0;
116         ConD2..           -Sum(i, b(i)  *(y1(i)+y2(i))) + Beta*q - Alpha*(w1+w2) =E= 0;
117         ConD3..                                        q - ThetaStar*(w1+w2) =E= 0;
118
119 Models  Primal_SCSC  / ObjP , ConP1, ConP2, ConP3 /
120         Dual_SCSC    / ObjD , ConD1, ConD2, ConD3 / ;
121
122 Solve Primal_SCSC using LP Maximizing Theta;
123       XbarStar(j) = (xbar1.L(j)+xbar2.L(j))/(w1.L+w2.L);
124       tStar       = p.L/(w1.L+w2.L);
125       UbarStar(i) = (ubar1.L(i)+ubar2.L(i))/(w1.L+w2.L);
126
127 Solve Dual_SCSC using LP Maximizing Theta;
128       YStar(i) = (y1.L(i)+y2.L(i))/(w1.L+w2.L);
129       zStar    = q.L/(w1.L+w2.L);
130       VStar(j) = (v1.L(j)+v2.L(j))/(w1.L+w2.L);
131
132 Put / /'Finding a Strict Comp. Solution via Approach I';
133 Put /'-------------------------------------------'/;
134 Put '       Primal              Dual     '/;
135 Put '-------------------------------------------'/;
136 Loop(j,
137       Put 'x(':>5; Put ord(j):<>3:0; Put ')= ':3; Put (XbarStar(j)/tStar):<10:3;
138       Put 'v(':>5; Put ord(j):<>3:0; Put ')= ':3; Put Vstar(j):<10:3;
139       Put /;
140     );
141 Put /;
142 Loop(i,
143       Put 'u(':>5; Put ord(i):<>3:0; Put ')= ':3; Put (UbarStar(i)/tStar):<10:3;
144       Put 'y(':>5; Put ord(i):<>3:0; Put ')= ':3; Put Ystar(i):<10:3;
145       Put /;
146     );
147 Put '-------------------------------------------'/ / /;
148
149 ***************
150 *First approach
151 **********************************************************
152
153 **********************************************************
154 *Second approach
155 ***************
```

```
156
157  Equations
158        ObjPD
159          ConPD ;
160
161        ObjPD..           Theta =E= Sum(j, xbar2(j)) + Sum(i, ubar2(i)) + Sum(i, y2(i)) + Sum(j,
v2(j)) + w2;
162          ConPD..                 Sum(j, c(j)*(xbar1(j)+xbar2(j))) + Alpha*p - q =E= 0;
163
164  Model PD_SCSC    / ObjPD, ConP1, ConP2, ConD1, ConD2, ConPD/ ;
165
166  Solve PD_SCSC using LP Maximizing Theta;
167         XbarStar(j) = (xbar1.L(j)+xbar2.L(j))/(w1.L+w2.L);
168         UbarStar(i) = (ubar1.L(i)+ubar2.L(i))/(w1.L+w2.L);
169         Ystar(i) = (y1.L(i)+y2.L(i))/(w1.L+w2.L);
170         Vstar(j) = (v1.L(j)+v2.L(j))/(w1.L+w2.L);
171
172  Put / /'Finding a Strict Comp. Solution via Approach I';
173  Put /'--------------------------------------------'/;
174  Put '         Primal                  Dual      '/;
175  Put '--------------------------------------------'/;
176  Loop(j,
177          Put 'x(':>5; Put ord(j):<>3:0; Put ')= ':3; Put (XbarStar(j)/tStar):<10:3;
178          Put 'v(':>5; Put ord(j):<>3:0; Put ')= ':3; Put Vstar(j):<10:3;
179          Put /;
180       );
181  Put /;
182  Loop(i,
183          Put 'u(':>5; Put ord(i):<>3:0; Put ')= ':3; Put (UbarStar(i)/tStar):<10:3;
184          Put 'y(':>5; Put ord(i):<>3:0; Put ')= ':3; Put Ystar(i):<10:3;
185          Put /;
186       );
187  Put '--------------------------------------------'/ / /;
188
189  ****************
190  *Second approach
191  ************************************************************
192
```